\documentclass[11pt]{article}
\usepackage[T1]{fontenc}
\usepackage{lmodern}
\usepackage[margin=1in]{geometry}
\usepackage{amsmath,amssymb,amsthm,mathtools}
\usepackage{microtype}
\usepackage{needspace}
\usepackage{enumitem}
\usepackage{xcolor}
\usepackage{booktabs}
\usepackage{aliascnt}
\usepackage[hidelinks]{hyperref}
\usepackage[nameinlink,noabbrev]{cleveref}
\usepackage{cite}
\makeatletter
\@ifundefined{dth@everypar}{}{%
  \def\cref@thmoptarg[#1]#2#3#4{%
    \cref@thmnoarg{#2}{#3}{#4}}}
\makeatother

\hypersetup{
  pdftitle={Weighted bisections of triangle-free subcubic graphs},
  pdfauthor={Shufei Wu and Qinghou Zeng},
  pdfkeywords={graph bisection; weighted cuts; triangle-free graphs; subcubic graphs; judicious partitions}
}

\newtheorem{theorem}{Theorem}[section]
\newaliascnt{lemma}{theorem}
\newtheorem{lemma}[lemma]{Lemma}
\aliascntresetthe{lemma}

\newaliascnt{conjecture}{theorem}
\newtheorem{conjecture}[conjecture]{Conjecture}
\aliascntresetthe{conjecture}

\newaliascnt{proposition}{theorem}
\newtheorem{proposition}[proposition]{Proposition}
\aliascntresetthe{proposition}

\newaliascnt{corollary}{theorem}

\aliascntresetthe{corollary}

\theoremstyle{definition}
\newaliascnt{definition}{theorem}

\aliascntresetthe{definition}

\theoremstyle{remark}
\newaliascnt{remark}{theorem}
\newtheorem{remark}[remark]{Remark}
\aliascntresetthe{remark}

\crefname{theorem}{Theorem}{Theorems}
\crefname{lemma}{Lemma}{Lemmas}
\crefname{conjecture}{Conjecture}{Conjectures}
\crefname{proposition}{Proposition}{Propositions}
\crefname{corollary}{Corollary}{Corollaries}
\crefname{definition}{Definition}{Definitions}
\crefname{remark}{Remark}{Remarks}
\crefname{section}{Section}{Sections}
\Crefname{section}{Section}{Sections}
\crefname{subsection}{Subsection}{Subsections}
\Crefname{subsection}{Subsection}{Subsections}

\newcommand{\claw}{K_{1,3}}
\newcommand{\Prb}{\mathbb P}

\newcommand{\R}{\mathbb R}
\newcommand{\U}{\mathcal U}

\newcommand{\conv}{\operatorname{conv}}

\setlist[enumerate]{itemsep=3pt,topsep=5pt}
\setlist[itemize]{itemsep=3pt,topsep=5pt}

\title{Weighted bisections of triangle-free subcubic graphs}

\author{
Shufei Wu$^1$\thanks{Research partially supported by
Natural Science Foundation of Henan Province
(Grant No. 252300421990).
Email: shufeiwu@hotmail.com}
\quad
Qinghou Zeng$^2$\thanks{Research partially supported by National Key R\&D Program of China (Grant No. 2023YFA1010202) and National Natural Science Foundation of China (Grant No. 12371342). Email: zengqh@fzu.edu.cn}\\
{\small $^1$School of Mathematics and Information Science, Henan Polytechnic University,}\\
{\small Jiaozuo 454003, Henan, China.}\\
{\small $^2$Center for Discrete Mathematics, Fuzhou University, Fuzhou 350003, Fujian, China}
}

\date{}

\begin{document}
\maketitle
\pagestyle{plain}

\begin{abstract}
We prove a conjecture of Gerke, Gutin, Yeo, and Zhou:
for every finite triangle-free subcubic graph other than the claw
and every nonnegative edge weighting, there is a bisection
whose weight is at least $11/15$ of the total edge weight. This bound is best
possible, as shown by the Petersen graph. The proof uses
matching-based local distributions, weighted induction, and
convex separation to establish the edgewise cut guarantees
while maintaining global balance.
\end{abstract}



\section{Introduction}
\label{sec:introduction}

For a finite graph $G$ and $X\subseteq V(G)$, let
$\delta_G(X)$ be the set of edges with exactly one endpoint in $X$.
This edge set is the \emph{cut} determined by the partition
$(X,V(G)\setminus X)$; its edges are said to be \emph{cut}.
A \emph{bisection} is such a partition $(X,Y)$ for which
$\bigl||X|-|Y|\bigr|\leq1$.
For a nonnegative edge weighting $w:E(G)\to\R_{\geq0}$,
write $w(F)=\sum_{e\in F}w(e)$ for $F\subseteq E(G)$;
the weight of the bisection $(X,Y)$ is $w(\delta_G(X))$.
A graph is \emph{subcubic} if its maximum degree is at most
three, and \emph{cubic} if every vertex has degree three.
Unless stated otherwise, graphs are simple.

For maximum cuts, Hopkins and Staton \cite{HopkinsStaton1982}
proved a \(4/5\) edge bound for triangle-free cubic graphs.
Bondy and Locke \cite{BondyLocke1986} extended the bound
to subcubic graphs, and Xu and Yu \cite{XuYu2008}
characterized the connected equality cases.
Zhu \cite{Zhu2009} obtained a $17/21$ bound for
$2$-connected triangle-free subcubic graphs, with seven exceptions. For nonnegative edge weights,
Gutin and Yeo \cite{GY2023} proved an $8/11$ cut bound for
triangle-free subcubic graphs. Hausmann \cite{Hausmann2026}
subsequently improved the bound to $15/19$.

For bisections, Bollob\'as and Scott \cite{BollobasScott2004}
established the best-possible $2/3$ bound for cubic graphs.
Lee, Loh, and Sudakov \cite{LLS2013} developed
Edwards-type bounds \cite{Edwards1973} for bisections. 
More recently, Guo and Gutin \cite{GuoGutin2026}
established sharp unweighted bisection bounds
for graphs of bounded maximum degree.
Further results under degree, matching, or forbidden-subgraph
conditions can be found in
\cite{MaWu2026,LinZeng2021,RHZ2022,HWZ2024,HuXu2025,HouYan2020}.
For related judicious partition problems, see
\cite{BollobasScott1999,BollobasScott2002,Scott2005,
XuYu2014,FHY2018,HouWu2021,WuXiong2024}.

Our question concerns the maximum weight of a bisection
in a triangle-free subcubic graph. Two examples show the
limits of any universal bound. Every bisection of the claw $K_{1,3}$ cuts
exactly two of its three edges; hence the claw must be
excluded from a bound greater than $2/3$. The Petersen graph has fifteen edges but no bisection
cutting more than eleven.
Gerke, Gutin, Yeo, and Zhou \cite{GGYZ2024} proved that
every edge-weighted subcubic graph has a bisection of weight
at least $2/3$ of its total edge weight. For triangle-free
subcubic graphs other than the claw, they established the
improved bound $613/855$ \cite{GGYZ2026} and formulated
the following conjecture.

\begin{conjecture}[Gerke et al. \cite{GGYZ2026}]\label{conj:GGYZ}
Every finite triangle-free subcubic graph other than
$K_{1,3}$ admits, for every nonnegative edge weighting,
a bisection of weight at least $11/15$ of the total edge weight.
\end{conjecture}

Gerke et al.\ \cite{GGYZ2026} proved the
conjectured bound for bridgeless triangle-free cubic graphs.
To settle the conjecture, one must also handle degree-two
vertices, bridges, disconnected graphs, and claw components
without losing the balance condition. 

\begin{theorem}\label{thm:target}
Let $G$ be a finite triangle-free subcubic graph with
$G\not\cong\claw$. For every weighting
$w:E(G)\to\R_{\geq0}$, there is a bisection $(X,Y)$ of $G$
such that
\[
  w(\delta_G(X))\geq\frac{11}{15}w(E(G)).
\]
The constant $11/15$ is best possible.
\end{theorem}

Thus, within the class of triangle-free subcubic
graphs, excluding $K_{1,3}$ raises the best possible
universal weighted-bisection constant from $2/3$
to $11/15$. We prove the equivalent distributional
assertion that a single distribution on bisections
cuts every edge with probability at least $11/15$.
The reduction in \cref{sec:reduction} leaves a
connected-graph theorem and a claw-absorption
theorem, proved using the matching-based
constructions developed below.

\paragraph{Notation and conventions.}
Throughout, set $\alpha=11/15$.
A \emph{bisection distribution} is a probability
distribution on the bisections of a graph.
For an odd-order graph, the larger part of a
bisection is called the \emph{majority side}.
We allow loopless multigraphs in intermediate constructions;
their edges are labelled and degrees count multiplicities.
A \emph{bundle} is the set of all labelled edges with
the same endpoints. For $S\subseteq V(H)$, write
$E_H(S)$ for the set of edges with both endpoints in $S$.
We use $\delta_H(S)$ as above, also for multigraphs.
Isolated vertices are included in balance calculations.

\section{Global reduction and the claw obstruction}
\label{sec:reduction}

The following equivalence allows us to prove
weighted bounds first and recover a common
bisection distribution afterward.

\begin{lemma}\label{lem:separation}
For a finite graph $G$ and $a\in\R$, the following are equivalent:
\begin{enumerate}[label=(\roman*)]
\item For every nonnegative edge weighting $w$, there is a
bisection of weight at least $a\,w(E(G))$.
\item There is a finite bisection distribution cutting each
edge with probability at least $a$.
\end{enumerate}
\end{lemma}
\begin{proof}
If $G$ has no edges, both assertions hold. Otherwise, let
$c_1,\ldots,c_N\in\{0,1\}^m$ be the incidence vectors of
all bisections of $G$, where $m=|E(G)|$. Set
\[
 P=\conv\{c_1,\ldots,c_N\},\qquad D=P-\R_{\geq0}^m.
\]
The set $D$ is convex and closed. To see closedness,
consider a convergent sequence $p_j-r_j$ with $p_j\in P$
and $r_j\in\R_{\geq0}^m$. By compactness of $P$, a
subsequence has $p_j\to p\in P$, and along that subsequence
$r_j$ converges to a nonnegative vector.

Suppose $b=a\mathbf1\notin D$. Choose a point
$y\in D$ nearest to $b$ and put $q=b-y\ne0$.
By the nearest-point property,
\[
q\cdot(z-y)\leq0\qquad(z\in D).
\]

Taking $z=y-te_j$ gives $q_j\geq0$ for every coordinate. Thus $q$ is an admissible nonnegative weighting, but
\[
 \max_i q\cdot c_i\leq q\cdot y<q\cdot b=a\sum_e q_e,
\]
contradicting the hypothesis, since
$q\cdot b-q\cdot y=\|q\|^2>0$.
Therefore $b\in D$, and some point of $P$ is coordinatewise
at least $b$. Its convex coefficients define the required
finite distribution.

Conversely, such a distribution gives expected weight at least $a\,w(E(G))$ for every $w\geq0$, so one of its outcomes attains that bound.
\end{proof}

The proof rests on two auxiliary results: a connected-graph
theorem and a claw-absorption theorem for a cubic component.

\begin{theorem}\label{thm:connected}
Let $G$ be a finite connected simple triangle-free
subcubic graph with $G\not\cong\claw$. For every
nonnegative real edge weighting $w$, $G$ has a bisection
of weight at least $\alpha w(E(G))$. Equivalently,
it has a finite bisection distribution cutting every
edge with probability at least $\alpha$.
\end{theorem}

\begin{theorem}\label{thm:claw-absorption}
For every connected simple triangle-free cubic graph $C$,
the disjoint union $C\sqcup\claw$ has an equal-part
bisection distribution cutting every edge with probability
at least $\alpha$.
\end{theorem}

\begin{proof}[Proof of \cref{thm:target}]
Assume first that $G$ has no claw component. Independently sample a bisection of each component from the distribution given by \cref{thm:connected}. Reverse whole components as necessary so that the majority sides of the odd-order components lie in the two global parts in numbers differing by at most one. Even-order components contribute no imbalance. Thus the resulting partition is a bisection of $G$, and the reversals leave all edge-cut indicators unchanged.

If $G$ has at least two claw components, connect them
in a path using auxiliary zero-weight edges between their leaves.
At each intermediate claw, use different leaves for its two connections. The resulting graph on those vertices is connected, simple, triangle-free, subcubic, and not a claw. Apply \cref{thm:connected} and delete the auxiliary edges. This gives an equal-part distribution on the union of all claw
components, which can then be combined with the remaining
components as above.

Suppose there is exactly one claw component. If another component $H$ is not cubic, choose a vertex of $H$ of degree at most two and join it to a claw leaf by an auxiliary edge. The graph on this union is a connected simple triangle-free subcubic graph of at least five vertices. Apply \cref{thm:connected} and delete the new edge. This includes the case that $H$ is an isolated vertex. Regard the resulting union as a single bisection block and combine it with the remaining components.

It remains to consider one claw together with at least one cubic component, since $G$ itself is not a claw. Apply \cref{thm:claw-absorption} to one cubic component together with the claw. All remaining cubic components have even order, so their distributions from \cref{thm:connected} combine directly with this equal-part block. This gives a common bisection distribution on $G$ cutting every edge with probability at least $\alpha$.

For any nonnegative weighting, the expected cut weight under the resulting common distribution is at least $\alpha w(E(G))$. At least one outcome has cut weight at least $\alpha w(E(G))$. Sharpness follows from \cref{prop:sharp}.
\end{proof}

\section{Distributional matching blocks}
\label{sec:blocks}

We construct balanced random partitions by splitting the
pairs of a matching, then merging selected pairs into blocks.

We use the following form of Edmonds' matching-polytope theorem \cite{Edmonds1965}. A vector $x\in\R^{E(H)}$ belongs to the convex hull of the matching incidence vectors of a simple graph $H$ if and only if
\begin{align}
x_e &\ge 0
    && (e\in E(H)), \notag\\
x(\delta_H(v)) &\le 1
    && (v\in V(H)), \notag\\
x(E_H(S)) &\le \frac{|S|-1}{2}
    && (S\subseteq V(H),\ |S|\text{ odd}).
    \label{eq:mp-odd}
\end{align}
A feasible vector can therefore be realized by a single
finite matching distribution with the prescribed edge marginals.
For applications of the matching-polytope method
to regular graphs, see Naddef and Pulleyblank
\cite{NP1981}.

\begin{lemma}\label{lem:multimatching}
The matching-polytope description above remains valid
for loopless multigraphs with labelled edges.
 Every bridgeless cubic multigraph has a finite distribution on perfect matchings in which each labelled edge has probability $1/3$.
\end{lemma}
\begin{proof}
For the first assertion, aggregate $x$ over each parallel bundle. The degree and odd-set sums are unchanged, so the aggregate vector has a matching decomposition in the simple underlying graph. When the matching decomposition selects an underlying edge
whose bundle has total mass $y>0$, choose its label $e$
with conditional probability $x_e/y$. Bundles of zero mass
are never selected. This produces the required labelled-edge marginals. Necessity follows by checking each matching.

For the second assertion, put $x_e=1/3$. Every vertex constraint is an equality. For an odd vertex set $S$ in a cubic graph,
\[
 3|S|=2|E_H(S)|+|\delta_H(S)|.
\]
Its boundary has odd size. It cannot have size one in a bridgeless graph, so its size is at least three. This gives \eqref{eq:mp-odd}. In a matching decomposition, each vertex is covered with probability one; hence every positive-probability matching is perfect. If two vertices are joined by $m$ labels, the probability that they form a matching pair is $m/3$, because the $m$ label-selection events are disjoint.
\end{proof}

For a loopless multigraph $H$ and a matching
$N$ covering all but at most one vertex,
let $Q=Q(H,N)$ be the simple graph obtained
by contracting each matching pair, retaining
the unmatched vertex if present, and
suppressing loops and parallel edges.
We call $Q$ the \emph{simple matching quotient}
of $(H,N)$, and the matching pairs and the
possible unmatched singleton the
\emph{initial blocks}.
When $H$ is subcubic, every vertex of $Q$
corresponding to a matching pair has degree at most four.

\begin{lemma}\label{lem:blocks}
Let $H$ be a loopless triangle-free
subcubic multigraph, and let $N$ be
a matching covering all but at most
one vertex. Write $Q=Q(H,N)$,
and let $K$ be a random matching of $Q$.

Then $K$ can be coupled with a random
bisection of $H$ such that every matching
pair is separated. If $N$ leaves a vertex
$s$ unmatched, then $s$ always lies on
the majority side. For an edge $e$
joining distinct initial blocks, with
image $h$ in $Q$,
\begin{equation}\label{eq:block}
\Prb(e\text{ is cut}\mid N)
=\frac12+\frac12\Prb(h\in K).
\end{equation}
All parallel labels with endpoints in
the same matching pair are cut with
probability one.
\end{lemma}

\begin{proof}
Let $P=\{a,a'\}$ and $P'=\{b,b'\}$ be two distinct matching pairs. Since $aa',bb'\in E(H)$ and $H$ is triangle-free, no endpoint of one pair is adjacent to both endpoints of the other. Thus, after parallel labels are ignored, the edges between $P$ and $P'$ form a matching of size at most two.

If there are two such edges, then after relabeling they are $ab$ and $a'b'$. Since each matching pair is split between the two sides, the condition that $ab$ be cut is equivalent to the condition that $a'b'$ be cut. Hence all labels between the two pairs impose the same requirement on their relative orientation. Thus a selected edge of $K$ allows its two pairs to be merged into a four-vertex block with a $2$--$2$ bipartition cutting every edge between the pairs.

A selected edge between the singleton $s$ and a pair gives a three-vertex block whose underlying simple graph is a path. Orient this block so that $s$ lies on its two-vertex side.
Unselected pairs remain two-vertex blocks, and an unselected singleton remains a one-vertex block. Since $K$ is a matching, no chain of mergers involves more than two initial blocks.

Orient each even block independently by a fair reversal. When an odd block exists, orient its larger side towards the global majority side. There is exactly one odd block, of order one or three, so every outcome is a bisection with the stated majority property.
If $h\in K$, every edge represented by $h$ is cut.
 Otherwise its endpoints are in distinct final blocks, at least one of which is an even block with an independent fair orientation. Its conditional cut probability is $1/2$, proving \eqref{eq:block}.
\end{proof}

The same block-reversal argument extends
to larger blocks of the simple matching quotient.
For a loopless triangle-free multigraph $H$
with a perfect matching $N$, let $Q=Q(H,N)$.
A set $B\subseteq V(Q)$ is \emph{compatible} if the union
of the matching pairs indexed by \(B\) admits a bipartition
that splits every pair and cuts every edge of $H$ joining
distinct pairs indexed by $B$.
A partition of $V(Q)$ is \emph{compatible} if each of its
blocks is compatible.

\begin{lemma}\label{lem:compatible-blocks}
Let $H$ be a loopless triangle-free multigraph with a perfect
matching $N$, and let $Q$ be its simple matching quotient.
Let $\Pi$ be a finite distribution on compatible partitions
of $V(Q)$. Then $\Pi$ lifts to a distribution on equal-part bisections of $H$. Every edge inside a matching pair is cut, and for every edge $e$ joining two distinct pairs with quotient endpoints $x,y$,
\begin{equation}\label{eq:compatible-block}
 \Prb(e\text{ is cut}\mid N)
 =\frac12+\frac12\Prb_{\Pi}(x,y\text{ lie in the same block}).
\end{equation}
In particular, every block inducing a tree in $Q$ is compatible.
\end{lemma}
\begin{proof}
After parallel labels are suppressed, the edges between two matching pairs form a matching of size at most two, by triangle-freeness. If two distinct edges occur, they impose the same requirement on the relative orientations of the pairs, as in the proof of \cref{lem:blocks}. Thus a block is compatible precisely when its relative-orientation constraints are consistent.

For each sampled partition and each of its blocks, choose one compatible bipartition and then reverse its two sides independently with probability $1/2$. Every matching pair remains split, so every outcome is an equal-part bisection and every edge inside a pair is cut. Let $e$ join distinct pairs with quotient endpoints $x,y$. If $x,y$ lie in the same block, the edge $e$ is cut by compatibility. If they lie in different blocks, the independent reversals of the two blocks make $e$ cut with probability $1/2$. This proves \eqref{eq:compatible-block}. If $Q[B]$ is a tree, its orientation requirements can be satisfied successively from a root, so $B$ is compatible.
\end{proof}

\begin{lemma}\label{lem:vectors}
Let $Q$ be simple with maximum degree at most four. Each of the following vectors belongs to its matching polytope:
\begin{enumerate}[label=(\roman*)]
\item $z_e=1/5$ for every edge;
\item $z_h=2/5$, all other entries $1/5$, if one endpoint of $h$ has degree at most three;
\item $z_e=2/5$ for all edges incident with a vertex $s$ of degree at most two, and $1/5$ elsewhere;
\item $z_h=3/5$, all other entries $1/5$, if both endpoints of $h$ have degree at most three.
\end{enumerate}
\end{lemma}
\begin{proof}
The degree inequalities follow directly: in (ii) a degree-four endpoint has total at most $1$, and in (iii) simplicity ensures that a neighbor of $s$ receives at most one increment. In (iv), each endpoint of $h$ has total at most $3/5+2/5=1$.

We check every odd-set inequality. Singletons have no internal edges. Write $k=|S|$. For (i), a three-set has at most three edges. For $k\geq5$, the degree bound gives $e_Q(S)\leq2k$, and
\[
 \frac{2k}{5}\leq\frac{k-1}{2}.
\]

For (ii), only sets containing both endpoints of $h$ need further attention. If $k=3$, the internal weight is at most $4/5$. If $k=5$, the low-degree endpoint sends at most three edges into the other four vertices, which span at most six edges, so the internal weight is at most $(9+1)/5=2$. For $k\geq7$ it is at most $(2k+1)/5\leq(k-1)/2$.

For (iii), consider sets containing $s$. At $k=3$ there are at most three internal edges and two increments, giving total at most $1$. At $k=5$ the other four vertices span at most six edges and $s$ contributes at most two, giving $(8+2)/5=2$. For $k\geq7$, the internal degree sum is at most $4(k-1)+2$, so at most $2k-1$ edges and two increments give the same upper bound $(2k+1)/5$.

For (iv), suppose $S$ contains both endpoints of $h$. At $k=3$ the bound is $(3+2)/5=1$. At $k=5$, the edge $h$, at most four edges from its endpoints to the remaining three vertices, and at most three edges within those vertices give at most eight edges; the total is at most $(8+2)/5=2$. For $k\geq7$ the two low-degree endpoints yield at most $2k-1$ internal edges, and the two increments again give $(2k+1)/5$.
All remaining odd sets contain no edge with an increased coordinate.
\end{proof}

\section{Interface-controlled local constructions}
\label{sec:carriers}

The $11/15$ bound for weighted bridgeless
triangle-free cubic graphs was proved by
Gerke et al.\ \cite{GGYZ2026}.
We next establish a distributional bound for
labelled cubic multigraphs, including the exact
cut probability of each labelled edge.

\begin{lemma}\label{lem:cubic}
Every loopless bridgeless triangle-free cubic
multigraph admits a bisection distribution in which
each labelled edge in a bundle of size $m$ is cut
with probability

\begin{equation}\label{eq:bundle-rate}
 \frac{m}{3}+\left(1-\frac{m}{3}\right)\frac35
 =\frac{9+2m}{15}\geq\alpha.
\end{equation}
\end{lemma}
\begin{proof}
Choose a perfect matching $M$ from \cref{lem:multimatching}. Conditional on $M$, use the constant $1/5$ vector of \cref{lem:vectors} on its simple quotient and apply \cref{lem:blocks}. Every edge within a matching pair is cut with probability one, while every edge joining distinct matching pairs is cut with probability $3/5$. The endpoints of a bundle of size $m$ are paired with probability $m/3$, giving \eqref{eq:bundle-rate}.
\end{proof}

The next lemma controls the probability that two specified
degree-two vertices lie on opposite sides, which we call their
\emph{terminal-separation probability}. We require bridgelessness
only after the two terminals are joined.

\begin{lemma}\label{lem:terminals}
Let $H$ be a connected loopless triangle-free multigraph with exactly two degree-two vertices $s,t$, each having two distinct neighbors, and with all other vertices of degree three. Suppose that $B=H+f$, where $f=st$ is a new labelled edge, is bridgeless.

There is a bisection distribution on $H$ cutting every edge with probability at least $\alpha$. If $s,t$ are adjacent or have no common neighbor, the same distribution satisfies
\[
 \Prb(s,t\text{ are separated})\geq\alpha.
\]
If they are nonadjacent and have a common neighbor, it instead satisfies
\[
 \Prb(s,t\text{ are separated})=\frac{8}{15}.
\]
\end{lemma}
\begin{proof}
If $s,t$ are adjacent or have no common neighbor, $B$ is triangle-free. Apply \cref{lem:cubic} to $B$ and delete the auxiliary label $f$, leaving the vertex partition unchanged. The marginal of $f$ gives the asserted terminal probability.

Now assume $s,t$ are nonadjacent and put
\[
 A=N_H(s)\cap N_H(t),\qquad F=\{sa,ta:a\in A\}.
\]
Here $1\leq|A|\leq2$, and all edges of $F$ are single labels. Every triangle of \(B\) has vertex set $W=\{s,t,a\}$ for some $a\in A$. Since $B$ is cubic and $W$ spans three edges, $3|W|=2|E_B(W)|+|\delta_B(W)|$ gives $|\delta_B(W)|=3$. Choose a perfect matching $M$ in $B$ with each label having marginal $1/3$. For such a set $W$, every perfect matching satisfies
\[
 2|M\cap E_B(W)|+|M\cap\delta_B(W)|=|W|=3.
\]
Hence $|M\cap\delta_B(W)|$ is a positive odd integer. Since $|\delta_B(W)|=3$ and each boundary label has matching marginal $1/3$,
\[
 \mathbb E|M\cap\delta_B(W)|=1.
\]
It follows that $|M\cap\delta_B(W)|=1$ for every matching in the support, and therefore $|M\cap E_B(W)|=1$. Since $E_B(W)=\{f,sa,ta\}$, every matching of positive probability simultaneously satisfies
\begin{equation}\label{eq:tight-triangle}
 |M\cap\{f,sa,ta\}|=1\qquad(a\in A).
\end{equation}

If $f\in M$, remove $F$. On the quotient of the triangle-free subcubic graph $B-F$ determined by $M$, use the constant $1/5$ matching vector from \cref{lem:vectors}(i) and then apply \cref{lem:blocks}. The pair $\{s,t\}$ and the pair containing any $a\in A$ have no edge between them: writing the latter pair as $\{a,c\}$, either $sc$ or $tc$ would create a triangle in $H$. They cannot be merged, even indirectly, because the quotient edges chosen for merging form a matching. Their final even blocks are independently oriented. Each deleted edge of $F$, evaluated on the resulting partition, therefore has cut probability $1/2$. Every ordinary edge has probability one if its endpoints are paired, and $3/5$ otherwise.

If $f\notin M$, then $M$ is a perfect matching of $H$. All edges of $F\setminus M$ map to one quotient edge $h$ with an endpoint of degree at most three. To see this when $A=\{a\}$, write $sa,ty\in M$ after exchanging $s,t$ if necessary. 
The sole edge of $F\setminus M$ is $at$.
There is no $sy$, which would make $y$ a second common neighbor, and no $ay$, which would form the triangle $aty$. The pair $\{s,a\}$ has at most three outgoing labels.
When $A=\{a,b\}$, \eqref{eq:tight-triangle} forces $M$ to use either $sa,tb$ or $ta,sb$. 
The two edges of $F\setminus M$ project to
 the same quotient edge $h$ and impose the same requirement on the relative orientations of the two pairs. Again an endpoint pair has at most three outgoing labels.

Give $h$ marginal $2/5$ and every other quotient edge marginal $1/5$, using \cref{lem:vectors}(ii). The edges of $F\setminus M$ have cut probability $7/10$, and all other interpair edges have probability $3/5$.

For a fixed $e\in F$, the disjoint events $e\in M$, $f\in M$, and neither have probabilities $1/3$ each. Consequently
\[
 \Prb(e\text{ is cut})=\frac13\left(1+\frac12+\frac7{10}\right)=\alpha.
\]
For an ordinary bundle of size $m$, the conditional probabilities remain one when paired and $3/5$ otherwise, giving \eqref{eq:bundle-rate}. Every outcome consists of balanced even blocks.

Finally, $f\in M$ forces $s,t$ to be separated. If $f\notin M$, selecting $h$ forces $s,t$ to have the same color, since every edge of $F$ is then cut. If $h$ is not selected, the two endpoint pairs lie in distinct final blocks, and $s,t$ are separated with probability $1/2$. Hence the terminal probability is
\[
 \frac13+\frac23\left(1-\frac25\right)\frac12=\frac8{15}.
\]

\end{proof}

For a unique degree-two vertex, Gerke et al.\
\cite{GGYZ2026} obtained bounds on
the probability that the vertex lies on the larger
side of a random bisection.
Here we need a distribution with the required edge-cut
marginals, exact root-edge probabilities, and a joint
condition involving the root and its neighbors.

\begin{lemma}\label{lem:rooted}
Let $H$ be a connected loopless bridgeless triangle-free multigraph with a unique degree-two vertex $r$, all other vertices having degree three. Then $H$ has odd order and a bisection distribution cutting every edge with probability at least $\alpha$. If $u,v$ are the neighbors of $r$, this same distribution satisfies
\begin{align}
\Prb(ru\text{ is cut})=\Prb(rv\text{ is cut})
  &=\frac34,\notag\\
\Prb(r\text{ is on the majority side})
  &=\frac{23}{30},\notag\\
r\text{ on the minority side}
  &\Longrightarrow
  u,v\text{ both on the majority side}.
  \label{eq:root-joint}
\end{align}
Every edge not incident with $r$, in a bundle of size $m$, has the cut probability given by \eqref{eq:bundle-rate}.
\end{lemma}
\begin{proof}
The degree-sum formula implies that the number of degree-three vertices is even. The two neighbors of $r$ are distinct: a double edge from $r$ to a single neighbor either forms an entire two-vertex component or leaves a bridge at that neighbor, both impossible here. They are nonadjacent by triangle-freeness. Suppress $r$ and insert the new labelled edge $f=uv$. The resulting multigraph $B$ is loopless and cubic. It is also bridgeless. Indeed, every edge of $H$ lies on a cycle. A cycle avoiding $r$ is unchanged, while a cycle containing $r$ necessarily uses the subpath $u-r-v$, which may be replaced by the new edge $uv$. Thus every surviving edge of $B$ lies on a cycle, and the new edge $f$ itself lies on the image of any cycle of $H$ containing $u-r-v$. Although $B$ may have triangles through $f$, \cref{lem:multimatching} still applies.

Choose its perfect matching $M$ with label marginals $1/3$. If $f\notin M$, view $M$ as a matching of $H$ missing only $r$. Its two neighbors belong to distinct matching pairs. Apply \cref{lem:vectors}(iii) to the quotient singleton $r$ and then \cref{lem:blocks}. The root is always on the majority side and its two incident edges each have probability $7/10$ of being cut. Ordinary interpair edges have probability $3/5$.

If $f\in M$, choose one of $u,v$ uniformly at random. Suppose $u$ is chosen, and let
\[
N=(M-\{f\})\cup\{ru\}.
\]
Then $N$ misses precisely $v$.
The quotient edge $h$ between $v$ and $\{r,u\}$ represents only $rv$ and has a degree-at-most-three endpoint $v$. Give it marginal $2/5$, using \cref{lem:vectors}(ii). The edge $ru$ is cut with certainty and $rv$ with probability $7/10$. If $h$ is selected, the odd block is $u-r-v$, with majority side $\{u,v\}$. Otherwise the root lies in an even block with an independent fair orientation. Thus the conditional probability that $r$ is on the majority side is $(1-2/5)/2=3/10$. The other choice is symmetric, so each root edge has conditional cut probability $17/20$.

Since $\Prb(f\in M)=1/3$, averaging gives
\[
 \frac23\frac7{10}+\frac13\frac{17}{20}=\frac34,
 \qquad
 \frac23+\frac13\frac3{10}=\frac{23}{30}.
\]
The endpoints of each such bundle are paired with probability $m/3$, and their interpair cut probability is $3/5$ in every case, proving the remaining marginal assertion.

For \eqref{eq:root-joint}, a minority root occurs only when $f\in M$. If the chosen root pair is $\{r,u\}$, the unmatched vertex $v$ is always on the majority side. When $h$ is selected, the odd block is $u-r-v$ with majority side $\{u,v\}$; when $h$ is not selected, the pair $\{r,u\}$ is split, so a minority $r$ forces $u$ onto the majority side. Thus in either case a minority root forces both $u$ and $v$ onto the majority side.
\end{proof}

\section{Completion of the bridgeless subcubic case}
\label{sec:bridgeless}

We handle general bridgeless subcubic graphs by adding
auxiliary edges, applying the preceding distributions,
and then deleting the auxiliary edges. 

\begin{theorem}\label{thm:bridgeless}
Every finite loopless bridgeless triangle-free subcubic multigraph has a finite bisection distribution cutting each labelled edge with probability at least $\alpha$.
\end{theorem}
\begin{proof}
We first consider a connected graph $G$. A one-vertex graph has no edges. A nontrivial bridgeless graph on two vertices consists of two or three parallel edges, and its vertices can be separated. Assume henceforth that $n=|V(G)|\geq3$. Every vertex has degree two or three. A degree-two vertex has two distinct neighbors: otherwise its sole neighbor either completes a two-vertex component or has a unique edge leaving that component, which would be a bridge.

First suppose that $G$ has even order. Let $D$ be the set of degree-two vertices. The handshake identity implies that $|D|$ is even. The induced graph $G[D]$ is simple and has maximum degree at most two. If it contains a cycle, that cycle has no edges leaving it, so it is all of $G$. Since $n$ is even, its bipartition is a bisection cutting every edge.

Otherwise $G[D]$ is a union of paths, including isolated vertices. Choose a maximum matching on each path and duplicate each chosen edge with an auxiliary parallel label. All covered vertices now have degree three, and the remaining degree-two vertices form an independent set. The added labels form two-cycles with their original edges, so the graph remains bridgeless and triangle-free.

As long as there are two remaining degree-two vertices
with no common neighbor, join them by an auxiliary edge.
 Since the remaining vertices are independent, the new edge creates no parallel pair; the absence of a common neighbor prevents a triangle.
Since the original graph is connected, the two vertices are already joined by a path, so the new edge lies on a cycle.
The process leaves an independent even set $R$ of degree-two
vertices whose two-element neighbor sets are pairwise
intersecting. Every other vertex has degree three.

We show that $|R|$ is zero or two. Fix $s\in R$ with neighbors $a,b$. Each of $a,b$ can be adjacent to at most two further vertices of $R$, giving $|R|\leq5$. It remains to rule out four. Four pairwise intersecting two-element sets cannot have a common element here, since that vertex would have degree at least four. If a pairwise intersecting family of two-element sets has no common element, take two distinct sets $\{a,b\},\{a,c\}$ and a set avoiding $a$; the latter must be $\{b,c\}$. Every member of the family is then one of these three sets. All three types occur. With four members and the degree bounds, their multiplicities must be $2,1,1$.

Thus two of $a,b,c$ have all three incident edges going to $R$, while the third has two. Since every pair among $a,b,c$ has a common neighbor in $R$, triangle-freeness forbids edges among these three vertices. The remaining edge at the third vertex is therefore the unique edge leaving $R\cup\{a,b,c\}$, contradicting bridgelessness.
 Hence $|R|\in\{0,2\}$.

If $R=\varnothing$, the resulting multigraph is cubic,
and \cref{lem:cubic} applies. If $|R|=2$, its two
remaining degree-two vertices satisfy
\cref{lem:terminals}; adding their closing edge
preserves bridgelessness.
 Delete all auxiliary labels from the resulting distribution. The vertex partitions and all original edge-cut indicators remain unchanged. This proves the even-order case.

Suppose next that $G$ has odd order, so $|D|$ is odd. If $|D|=1$, use \cref{lem:rooted}. If $|D|\geq3$, some two vertices $a,b\in D$ are nonadjacent, since otherwise three of them form a triangle. Add a new vertex $x$ and edges $xa,xb$. The graph remains subcubic and triangle-free, and the new path lies on a cycle with an $a$--$b$ path in $G$. It is therefore bridgeless and has even order. Apply the even-order case and remove $x$ from each outcome. Removing one vertex from an equal partition gives a   bisection of the original graph. Every original cut indicator is preserved.

Finally, if $G$ is disconnected, apply the connected result independently to every component, including isolated vertices. Even components are already equally divided. Orient the odd components so that the numbers whose majority sides lie in the two global parts differ by at most one. The total imbalance is consequently zero or one. Whole-component reversals preserve all edge-cut indicators. The empty graph has its unique empty partition. This completes the proof.
\end{proof}

\begin{remark}
Parallel labels must be retained during completion: merging them can turn a bridgeless multigraph into a simple graph with bridges.
\end{remark}

\section{The single-leaf interface}
\label{sec:leaf}

The following lemma strengthens the
corresponding single-leaf bound of
Gerke et al.\ \cite{GGYZ2026}.

\begin{lemma}\label{lem:single-leaf}
Let $H$ be a connected simple bridgeless triangle-free graph whose degree-two vertices are precisely $v,a,b$, with $N_H(v)=\{a,b\}$, and whose other vertices have degree three. Add a new leaf $u$ adjacent to $v$, obtaining $G$. Then $G$ has an equal-part bisection distribution cutting every edge, including $uv$, with probability at least $\alpha$.
\end{lemma}
\begin{proof}
Write $N_H(a)=\{v,x\}$ and $N_H(b)=\{v,y\}$. By triangle-freeness, $x,y$ lie outside $\{a,v,b\}$.
They are distinct: otherwise the four-cycle $xavbx$ would
have a unique outgoing edge at $x$, which would be a bridge.
Suppress $a$ and $b$, and relabel $v$ as $r$.
The path $x-a-v-b-y$ becomes $x-r-y$.
 The resulting graph $L$ is connected, bridgeless, and simple,
with $r$ as its unique degree-two vertex. Its only possible
new triangle is $xry$, which occurs exactly when $xy\in E(H)$.

We use the following balanced four-vertex replacement in both cases. Given fixed colors on $x,y$, color
the four vertices $a,v,b,u$ with two vertices of each
color as follows. Uniformly choose one of $va,vb,vu$ to be the unique uncut star edge. There are two orientations of the resulting $2$--$2$ coloring; choose one that maximizes the number of cut edges among $xa,by$, breaking ties fairly. Each star edge is cut with probability $2/3$. If $x,y$ have the same color, each outer edge has cut probability $2/3$: the choice of $vu$ leaves both outer edges cut, while each of the other two choices gives probability $1/2$ to each. If $x,y$ have different colors, these probabilities are instead $5/6$: the choice of $vu$ gives each probability $1/2$, and the other two choices cut both. The replacement leaves all outside colors unchanged and
adds exactly two vertices to each part.

\medskip\noindent\emph{Case 1: $xy\notin E(H)$.}

Here $L$ is triangle-free. Apply \cref{lem:rooted}, orienting its majority side as $+$. Its joint property gives
\[
 \Prb(r=-)=\frac7{30},\qquad r=-\ \Longrightarrow\ x=y=+.
\]
Remove $r$ and insert $a,v,b,u$. If $r=-$, the outside imbalance is $+2$. Set $v=+$ and $a=b=u=-$. This contributes imbalance $-2$ and cuts all five local edges $xa,av,vb,by,uv$. If $r=+$, the outside is equally divided; use the balanced replacement just described. Thus each local edge is cut with probability at least
\[
 \frac7{30}+\frac{23}{30}\frac23=\frac{67}{90}>\alpha.
\]
All other cut indicators are unchanged. Every resulting partition of $G$ has equal parts.

\medskip\noindent\emph{Case 2: $xy\in E(H)$.}

Let $C=H-\{a,v,b\}$ and add a second label $f$ parallel to $xy$. This gives a loopless bridgeless triangle-free cubic multigraph $B$: it is obtained by suppressing the path $x-a-v-b-y$ to an additional $xy$ edge. Let $xp,yq$ be the other edges at $x,y$ in $B$. The vertices $p,q$ lie outside $\{x,a,v,b,y,u\}$ and are distinct, since $p=q$ would form a triangle with $xy$.

Choose a perfect matching $M$ of $B$ from \cref{lem:multimatching}.
Let $T$ be the event that $x,y$ form a matching pair.
The two parallel labels each belong to $M$ with probability $1/3$,
and these events are disjoint. Hence
\begin{equation}\label{eq:leaf-T}
 \Prb(T)=\frac23,\qquad \Prb(T^c)=\frac13.
\end{equation}

If $T^c$ occurs, the matching contains
$xp,yq$ and is a perfect matching of $C$. In its simple quotient, the two pairs $\{x,p\},\{y,q\}$ each have at most three outgoing labels, because $x,y$ have degree two in $C$. Their quotient edge $h$ represents $xy$ and possibly $pq$. No $xq$ or $yp$ is possible by triangle-freeness, so the represented requirements are compatible. Apply \cref{lem:vectors}(iv), assigning $h$ marginal $3/5$, and use pair blocks.
The resulting partition of $C$ has equal parts,
$xy$ is cut with probability $4/5$, and every
ordinary nonmatching edge has conditional cut
probability at least $3/5$.

Insert the four-vertex balanced replacement above. Each star edge has probability $2/3$, while each outer edge $xa,by$ has probability
\[
 \frac15\frac23+\frac45\frac56=\frac45.
\]
Every cut indicator in $C$ is unchanged.

If $T$ occurs, replace the matching
pair $\{x,y\}$ by the six-vertex block
\[
 S=\{x,a,v,b,y,u\}.
\]
Choose a balanced template from the following distribution. The displayed part has three vertices; the other part is its complement in $S$.
\begin{center}
\begin{tabular}{ccc}
\toprule
Probability & One part of the template & Uncut local edges\\
\midrule
$3/10$ & $\{x,v,y\}$ & $xy$\\
$7/30$ & $\{x,b,u\}$ & $va$\\
$7/30$ & $\{x,v,b\}$ & $vb$\\
$7/60$ & $\{y,v,u\}$ & $uv,xa$\\
$7/60$ & $\{x,v,u\}$ & $uv,by$\\
\bottomrule
\end{tabular}
\end{center}
The table gives
\begin{equation}\label{eq:leaf-template}
\begin{aligned}
\Prb(va\text{ is cut})
=\Prb(vb\text{ is cut})
=\Prb(uv\text{ is cut})
&=\frac{23}{30},\\
\Prb(xy\text{ is cut})
&=\frac7{10},\\
\Prb(xa\text{ is cut})
=\Prb(by\text{ is cut})
&=\frac{53}{60}.
\end{aligned}
\end{equation}

These indicators are unchanged by reversing the entire block.

Use $S$ and the remaining matching pairs as initial blocks. Their simple quotient $Q$ has degree at most two at $S$ and at most four elsewhere. A quotient matching merges initial blocks only in pairs.
In the conditional matching distributions specified
below, every selected quotient edge represents a
bundle with compatible orientation requirements.
Orient its two incident initial blocks so that
every edge in that bundle is cut.
  Then reverse each merged or unmerged block independently
and fairly. Every initial block is balanced, so all outcomes are balanced. A bundle not selected joins distinct final blocks and its edges have probability $1/2$ of being cut. Ordinary pair-to-pair bundles are compatible by triangle-freeness.

If $p,q$ belong to different matching pairs, the two edges leaving $S$ lead to different blocks and have no compatibility obstruction. Use the constant $1/5$ quotient vector, conditional on each template. Both outgoing edges $xp,yq$ have probability $3/5$, as do all other ordinary nonmatching edges.

If $p,q$ belong to the same matching pair, the sole quotient edge $h$ at $S$ represents both $xp$ and $yq$. Because $p,q$ have opposite colors, both requirements are compatible exactly when $x,y$ have opposite colors. This template event has probability $7/10$. Conditional on it, assign $h$ marginal $2/7$ and every other quotient edge marginal $1/5$. This follows by convexity from
\cref{lem:vectors}(i),(ii), since $h$ has a degree-one
endpoint and $1/5\leq2/7\leq2/5$.
On the other templates, assign $h$ marginal zero and
all other entries $1/5$, using the downward-closedness
of the matching polytope. Thus an incompatible bundle is never selected. Each outgoing edge has unconditional cut probability
\[
 \frac7{10}\left(\frac12+\frac12\frac27\right)
 +\frac3{10}\frac12=\frac35.
\]
All other quotient marginals remain $1/5$ even conditional on the template. This specifies a single joint law: sample $M$, then the
six-vertex template and its conditional quotient matching,
and finally the independent block reversals.

By \eqref{eq:leaf-T} and \eqref{eq:leaf-template}, the cut probabilities of a star edge and of \(xy\), respectively, are
\[
 \frac23\frac{23}{30}+\frac13\frac23=\alpha,
 \qquad
 \frac23\frac7{10}+\frac13\frac45=\alpha.
\]
Each outer edge $xa,by$ has probability
\[
 \frac23\frac{53}{60}+\frac13\frac45=\frac{77}{90}.
\]
Every other original edge is a single label in $B$ and belongs to $M$ with probability $1/3$. If its endpoints are paired by $M$, it is certainly cut; otherwise, its cut probability is at least $3/5$. This includes the case where $p,q$ form a matching pair.
Hence its total probability is at least $1/3+(2/3)(3/5)=\alpha$.
\end{proof}

\section{Gluing across bridge blocks: the connected theorem}
\label{sec:connected}

We prove \cref{thm:connected} by weighted
induction across bridge blocks, following
the reduction of Gerke et al.\ \cite{GGYZ2026}
with the stronger interface estimates
established above.
Working with weights allows the induction to proceed directly;
\cref{lem:separation} then recovers a single bisection distribution.

A \emph{bridge block} is a connected component after all
bridges are deleted; it may consist of one vertex.
Contracting the blocks gives a tree, since a cycle in the
contraction would place a bridge on a cycle. Write $c(G)$
for the number of blocks, one more than the number of
bridges. Every nontrivial block is bridgeless, with internal
degrees two or three. Such a block is \emph{tight} if it has
exactly one vertex of internal degree two. It then has odd
order, and in a connected graph with bridges can only be
a leaf block: its unique degree-two vertex is the only
possible bridge attachment.

\begin{proof}[Proof of \cref{thm:connected}]
We induct lexicographically on $(|V(G)|,c(G))$,
within the class of connected graphs other than the claw. The connected triangle-free graphs of order at most three have a bisection that cuts every edge. The bridgeless case follows from \cref{thm:bridgeless}. Assume that $|V(G)|\geq4$ and $c(G)\geq2$.

\medskip\noindent\emph{Case 1: an even leaf block.}

Suppose an even leaf block $H$ is attached by the bridge $ha$, with $h\in V(H)$, to the connected remainder $R$. If $R\not\cong\claw$, apply induction to $R$ and \cref{thm:bridgeless} to $H$. Reverse $H$ to cut $ha$. The resulting partition is a bisection, and its weight is at least
\[
 \alpha\bigl(w(E(H))+w(E(R))\bigr)+w(ha)
 \geq\alpha w(E(G)).
\]

If $R$ is a claw with center $z$ and leaves $a,b,d$, its bridge endpoint $a$ is a leaf because the center already has degree three. Within $H$, the number of degree-two vertices is even. Since $h$ is one of them, there is another one $t$. Only $h$ meets the bridge to $R$, so $t$ still has degree two in $G$. Add a zero-weight edge $tb$. Its endpoints have no common neighbor: all neighbors of $t$ are in $H$ and the sole neighbor of $b$ is $z$. The new edge completes a cycle through a $t$--$h$ path, $ha,az,zb$, so the bridge count decreases. The order is unchanged, and the graph contains a cycle and therefore is not a claw. Apply induction and delete the auxiliary edge.

We may now assume that all leaf blocks have odd order.

\medskip\noindent\emph{Case 2: no tight block.}

Assume that $G$ has no tight block.
Suppose there is a nontrivial leaf block $H$ with bridge root $h$. Because $H$ is not tight, it has an internal degree-two
vertex $t\ne h$, which is incident with no bridge. If the block tree has more than two vertices, choose a second leaf block $J$; it is not adjacent to $H$. If $J$ is nontrivial, choose a degree-two vertex $s$ other than its bridge root; if $J$ is a singleton, let $s$ be its vertex. In both cases $s,t$ have spare degree, are nonadjacent, and have no common neighbor. Indeed, all neighbors of $t$ are in $H$, and a common neighbor would yield an edge directly from $J$ to $H$. Add $st$ with weight zero and apply induction with fewer bridges.

The same argument works when the block tree has exactly two nontrivial blocks. Choose $s,t$ away from their unique connecting bridge. A common neighbor in either block would force one of $s,t$ to be an endpoint of that bridge, contrary to the choice.

The remaining configuration with a nontrivial leaf block
is $G=H+hy$, where $y$ is a new leaf. Let $D$ be the degree-two vertices of $H$. Its order is odd and it is not tight, so $|D|\geq3$ is odd. If some $t\in D$ lies outside $N_H[h]$, add the zero-weight edge $ty$. This creates no triangle, and it puts the only bridge $hy$ on a cycle. The new graph is bridgeless, so \cref{thm:bridgeless} applies. If no such $t$ exists, then $D\subseteq N_H[h]$.
Since $N_H[h]$ has three vertices and $|D|\geq3$, we have
\[
 D=N_H[h]=\{h,a,b\}.
\]
This is exactly \cref{lem:single-leaf}, which supplies the required bisection.

Finally, suppose every leaf block is a singleton. These are degree-one vertices of $G$. If two leaves have no common neighbor, join them by a zero-weight edge and use the reduction in bridge count. Otherwise every pair of leaves has a common neighbor. Since each leaf has exactly one neighbor, all leaves are adjacent to the same vertex $z$. There are at least two leaves. The vertex $z$ must itself be a singleton bridge block, since two leaf edges and the at least two internal edges of a nontrivial block would exceed degree three. If the block tree had another branch from $z$, following that branch to its end would produce a leaf block not adjacent to $z$. This is impossible. Hence $G$ consists of $z$ and its leaves; its order at least four and maximum degree at most three force $G=\claw$, the excluded case.

\medskip\noindent\emph{Case 3: a tight leaf block.}

Let $S$ be a tight leaf block with unique internal degree-two vertex $u$, joined by the bridge $uv$ to the connected remainder $T$. The block $S$ has odd order greater than three. Write $b=w(uv)$, $W_S=w(E(S))$, and $W_T=w(E(T))$. Sample a bisection of $S$ from the single joint distribution of \cref{lem:rooted}. Its expected cut weight is at least $\alpha W_S$, and
\[
 p:=\Prb(u\text{ is on the majority side})=\frac{23}{30}.
\]

If $T=\{v\}$, put $v$ on the minority side of $S$. This gives an equal-part partition and cuts the bridge with probability $p\geq\alpha$.

If $T$ is a claw, its bridge endpoint $v$ is a leaf; write $z$ for its center and $a,b$ for the other leaves. When $u$ is on the majority side of $S$, color the added vertices to cut all four edges $uv,vz,za,zb$. When $u$ is on the minority side, choose one of these four edges uniformly to be the unique uncut edge, and cut the other three. Each choice extends the sampled partition to a bisection of $G$: the all-cut extension places one new vertex on the majority side of $S$ and three on its minority side, and the four extensions with one uncut edge place either one
or two new vertices on its majority side. Hence every new edge is cut with probability $(3+p)/4\geq\alpha$.

If $Q:=G[V(T)\cup\{u\}]$ is a claw, then $T$ is the path with center $v$ and leaves $a,b$. When $u$ is on the minority side, cut all three new edges $uv,va,vb$. When $u$ is on the majority side, choose one of these edges uniformly as the unique uncut edge. Each of these extensions is balanced: precisely one of $v,a,b$ is added to the majority side of $S$. Every new edge is therefore cut with probability $1-p/3\geq\alpha$.
In all three cases the original edge marginals of $S$
are preserved, so the weighted bound follows.

Assume that $T$ is neither a single vertex nor a claw,
and that $Q=G[V(T)\cup\{u\}]$ is not a claw.
By induction, fix a bisection of $T$ of weight at least
$\alpha W_T$.

First suppose that $T$ has even order, or that $T$
has odd order and $v$ lies on its majority side.
Extend the sampled partition of $S$ to $P=G[V(S)\cup\{v\}]$ by placing $v$ on the minority side of $S$. The graph $P$ is bisected and its expected cut weight is at least $\alpha(W_S+b)$, since the bridge is cut with probability $p$. Identify the copies of $v$ in $P$ and $T$, reversing the coloring of $P$ if necessary. The union is balanced: its imbalance has absolute value one when $T$ is even, and zero when $T$ is odd because $v$ belongs to the majority side of $T$. The edge weights add, proving the bound.

It remains that $T$ has odd order and $v$ lies on its minority side in the chosen bisection. By induction, take a bisection of $Q$ cutting weight at least $\alpha(W_T+b)$. Sample the rooted distribution on $S$. If $u$ is on the majority side, identify the two copies of $u$ in $S$ and $Q$, orienting the bisection of $Q$ so that its color at $u$ agrees with the majority color of $S$. The union is balanced, because $Q$ has even order and the signed imbalance is $1+0-1=0$. If $u$ is on the minority side, combine $S$ with the fixed bisection of $T$, orienting their majority sides oppositely. This gives an equal-part partition and cuts the bridge, since $u$ and $v$ are both minority vertices of their respective odd-order pieces. The two extensions preserve the sampled cut inside $S$. Taking expectations, their total cut weight is at least
\[
 \alpha W_S+p\alpha(W_T+b)+(1-p)(\alpha W_T+b)
 =\alpha w(E(G))+(1-p)(1-\alpha)b\geq\alpha w(E(G)).
\]
This completes the weighted induction. The distributional
assertion now follows from \cref{lem:separation}.
\end{proof}

\section{Claw absorption and three-edge-connected cores}
\label{sec:absorption}

We now prove the claw-absorption theorem from
\cref{sec:reduction}. For a connected simple triangle-free
cubic graph $C$, write $\U(C)$
for the assertion of
\cref{thm:claw-absorption}. We first deal with cubic graphs
having a bridge, then reduce two-edge cuts to
three-edge-connected cores. 
It remains to construct a suitable distribution
for each such core.

\subsection{Cubic graphs with a bridge}

Sampling the two sides independently need not give the
required bridge-cut probability. We therefore couple their
root states.

\begin{lemma}\label{lem:bridged-absorption}
Let $C$ be a connected simple triangle-free cubic graph. If $C$ has a bridge, then $\U(C)$ holds.
\end{lemma}
\begin{proof}
Delete a bridge $xy$ and let $H_1,H_2$ be its two connected sides, rooted at $x,y$ respectively. Each root has degree two in its side and every other vertex has degree three, so each side has odd order. Add a new leaf to each root, obtaining connected even-order graphs $F_1,F_2$. They are simple triangle-free subcubic graphs, and neither is a claw: a three-vertex $H_i$ could not have the stated degree sequence, while a single root could not have degree two.

Apply \cref{thm:connected} separately to $F_1,F_2$ and delete their new leaves. On $H_i$ this gives a bisection distribution. Let $R_i$ be the indicator that its root is on its majority side.
Deleting the new leaf makes the color class opposite to that leaf the majority side of $H_i$. Hence $R_i=1$ exactly when the deleted leaf edge was cut. Thus
\[
 p_i:=\Prb(R_i=1)\geq\alpha,
\]
in the same distribution that retains all internal edge marginals.

Couple $R_1$ and $R_2$ to maximize their agreement,
then sample each side conditionally on its root state.
 The coupling preserves both side distributions and gives
the disagreement event $D$ probability
\[
 \delta=|p_1-p_2|\leq1-\alpha=\frac4{15}.
\]
Choose an event $A$ of probability $t=1/5$ having maximum overlap with $D$: use $D\subseteq A$ if $\delta\leq t$, and $A\subseteq D$ otherwise. Splitting finitely many atoms realizes this choice without changing either sampled partition or either side's marginal distribution; it only refines the probability space. Then
\[
 \Prb(A\mathbin{\triangle}D)=|t-\delta|\leq\frac15.
\]

On $A$, orient the majority sides of $H_1,H_2$ in the same direction, giving $C$ imbalance two. Cut every edge of the claw, orienting its three leaves into the minority side of $C$; the full graph is equally divided. Outside $A$, orient the two majority sides in opposite directions and choose a uniformly random $2$--$2$ partition of the claw, of which there are three up to reversal. Again the full graph is equally divided.

If the root states agree, the bridge is cut exactly when the two majority sides point in opposite directions. If they disagree, the opposite criterion holds. Thus the bridge is uncut exactly on $A\mathbin{\triangle}D$, and has cut probability at least $4/5$. Every claw edge has probability
\[
 \frac15+\frac45\frac23=\alpha.
\]
Whole-side reversals preserve all internal cut indicators. This proves $\U(C)$.
\end{proof}

\subsection{Reduction at a two-edge cut}

The joining-edge cut probabilities depend only on the separation states of the two terminal pairs, not on their full internal colorings.

\begin{lemma}\label{lem:parity-coupling}
Let $H_1,H_2$ be vertex-disjoint graphs, each with a distribution on equal-part bisections and with an ordered terminal pair $(s_i,t_i)$ in $H_i$. Let $p_i$ be the probability that the terminals of $H_i$ are separated. The two distributions can be coupled and oriented so that, after adding the edges $s_1s_2$ and $t_1t_2$, both joining edges have cut probability $1-\frac12|p_1-p_2|$. Every original edge retains its cut distribution, and the union remains equally divided.
\end{lemma}
\begin{proof}
Couple the two binary terminal parities with maximum agreement, and sample each block conditionally on its parity. The parity disagreement probability is $|p_1-p_2|$. If they agree, a relative reversal cuts both joining edges. If they disagree, the two possible relative orientations each cut one joining edge, so choose them fairly. All blocks have equal part sizes, so these reversals preserve balance. A final fair global reversal may be used if fully labelled orientations are wanted.
\end{proof}

\begin{lemma}\label{lem:minimal-two-shore}
Let $C$ be a connected simple bridgeless triangle-free cubic graph with a two-edge cut. Among all nonempty proper $S\subset V(C)$ with $|\delta_C(S)|=2$, choose one of minimum order. Let $H=C[S]$, with terminals $s,t$ incident with the cut edges. Then $H$ is connected, $s\ne t$, $|S|$ is even and at least six, and $K=H+st$ is a simple three-edge-connected cubic graph. Moreover, $K$ is triangle-free if and only if $s,t$ have no common neighbor in $H$.
\end{lemma}
\begin{proof}
Deleting the two cut edges leaves exactly two connected sides: every resulting component has at least two edges to its complement because $C$ has no bridge. Both cut edges cannot meet the same vertex $s$ of $S$, since its sole remaining edge would be a bridge; a singleton side would have three cut edges. Thus $H$ is connected with distinct degree-two terminals $s,t$, all other vertices having degree three. The degree sum $3|S|-2$ shows that $|S|$ is even. It cannot be two or four: two vertices do not have these degrees in a simple graph, and four vertices would require five edges, impossible in a triangle-free graph on four vertices. Hence $|S|\ge6$.

Fix a nonempty proper $T\subset S$. Replace $T$ by $S\setminus T$ if necessary so that $T$ does not contain both terminals.
Then $|T\cap\{s,t\}|\le1$. The new edge $st$
belongs to $\delta_K(T)$ exactly when $T$
contains one terminal. Each terminal in $T$
also contributes one cut edge to $\delta_C(T)$.
Hence
\[
|\delta_K(T)|
=|\delta_H(T)|+|T\cap\{s,t\}|
=|\delta_C(T)|.
\]
The last boundary is nonzero since $C$ is connected, cannot have size one because $C$ has no bridge, and cannot have size two because $0<|T|<|S|$ contradicts minimality. Therefore $K$ is three-edge-connected.

If $st$ were already an edge of $H$, its two labels in $K$ would make $|\delta_K(\{s,t\})|=2$, contradicting three-edge-connectivity since $|S|\ge6$. Thus $K$ is simple. As $H$ is triangle-free, any triangle of $K$ must contain the added edge $st$; it exists precisely when $s,t$ have a common neighbor in $H$.
\end{proof}

\begin{lemma}\label{lem:two-edge-reduction}
Let $C$ be a connected simple bridgeless triangle-free cubic graph with a two-edge cut. Assume that $\U(K)$ holds for every simple triangle-free three-edge-connected cubic graph $K$. Then $\U(C)$ holds, without any induction on $|V(C)|$.
\end{lemma}
\begin{proof}
Choose a minimum-order shore $H$ as in \cref{lem:minimal-two-shore}, with terminals $s,t$, and write $J$ for the opposite shore, with terminals $s',t'$. The joining edges are $ss',tt'$. Both shores have even order. We construct a balanced distribution on each side, placing the extra claw on one side, and couple their terminal parities by \cref{lem:parity-coupling}. The joint distribution within each side is preserved.

If $s,t$ have a common neighbor, \cref{lem:terminals} gives a balanced distribution on $H$ with all internal edge marginals at least $\alpha$ and terminal separation probability $p=8/15$. Join a degree-two terminal of $J$ to a leaf of the extra claw by an auxiliary edge of weight zero. The resulting graph is connected, simple, triangle-free, subcubic, of even order, and not itself a claw. Apply \cref{thm:connected} and discard the auxiliary edge. This gives a balanced distribution on $J\sqcup\claw$ with all its original edges cut with probability at least $\alpha$; its terminal separation probability $q$ may be arbitrary. Since $|8/15-q|\le8/15$, \cref{lem:parity-coupling} cuts each joining edge with probability at least $1-4/15=\alpha$.

Otherwise, \cref{lem:minimal-two-shore} gives a simple triangle-free three-edge-connected cubic graph $K=H+st$.
By hypothesis, $\U(K)$ holds, so
$K\sqcup\claw$ admits an equal-part
bisection distribution cutting every edge
with probability at least $\alpha$.
Deleting the auxiliary edge $st$ leaves a balanced
distribution on $H\sqcup\claw$ with all original
edge marginals at least $\alpha$ and terminal
separation probability $p\ge\alpha$.
 Meanwhile $J+s't'$ is a bridgeless cubic multigraph: every bridge of $J$ separates its terminals and lies on a cycle after the edge is added.
Thus \cref{lem:terminals} gives a balanced
distribution on $J$ with all internal edge marginals
at least $\alpha$ and terminal separation
probability $q\ge8/15$, whether its terminals are adjacent, share a common neighbor, or neither. Since $p\in[\alpha,1]$ and $q\in[8/15,1]$,
\[
 |p-q|\le\frac7{15},
\]
so the joining edges each have cut probability at least $1-7/30=23/30>\alpha$.

In both cases the coupling uses entire balanced block distributions and only reverses the relative color names. The union is a   bisection of $C\sqcup\claw$ with the required simultaneous edge marginals.
\end{proof}

\subsection{Three-edge-connected cores: the uniform construction}

We first show how to absorb a claw by occasionally
flipping one vertex of a balanced partition. We then
combine this operation with matching distributions on
three-edge-connected cubic graphs.

\begin{lemma}\label{lem:flip-absorption}
Let $C$ be a finite loopless graph of even order $n\geq2$, with a finite distribution on equal-part bisections. Write $p_e$ for the cut probability of each edge $e$ in this distribution. Then $C\sqcup\claw$ has one finite distribution on equal-part bisections in which each claw edge is cut with probability $\alpha$ and each edge $e$ of $C$ is cut with probability
\begin{equation}\label{eq:flip}
 p'_e=\left(1-\frac4{5n}\right)p_e+\frac2{5n}.
\end{equation}
The restriction to $C$ has equal parts with probability $4/5$ and part sizes $n/2-1,n/2+1$ with probability $1/5$. Parallel labelled edges are allowed.
\end{lemma}
\begin{proof}
Sample a bisection of $C$ from the given law and, independently, a mode $Z$ with $\Prb(Z=0)=4/5$ and $\Prb(Z=1)=1/5$. In mode zero retain the sampled bisection and independently choose a uniformly random balanced claw partition. Each claw edge is then cut with probability $2/3$.

In mode one choose a vertex of $C$ uniformly and independently of both the sampled bisection and $Z$, and flip its color. The two parts of $C$ now have sizes $n/2-1,n/2+1$. Put the claw center into the larger part and all three leaves into the smaller part. Each resulting part of the union has size $n/2+2$, and all three claw edges are cut. Mode zero is also equally divided. Thus every claw edge has total cut probability $(4/5)(2/3)+1/5=\alpha$.

Because $e$ has two distinct endpoints, the flip toggles its cut indicator with probability $2/n$, independently of that indicator. Hence
\[
 p'_e=\frac45p_e+\frac15\left[\left(1-\frac2n\right)p_e+\frac2n(1-p_e)\right],
\]
which is \eqref{eq:flip}. All choices have finite support and are shared by all edges, so one joint law realizes these marginals for every labelled edge.
\end{proof}

\begin{lemma}\label{lem:six-eight}
If $C$ is a connected simple bridgeless triangle-free cubic graph of order six or eight, then $\U(C)$ holds.
\end{lemma}
\begin{proof}
Choose a perfect matching from the distribution in
\cref{lem:multimatching}.
 Its simple quotient has $q=|V(C)|/2\in\{3,4\}$ vertices. The vector with every quotient entry $1/3$ is in the matching polytope: degrees are at most three, and the only nonsingleton odd sets have size three and contain at most three edges. Pair blocks give conditional probability $2/3$ for every nonmatching edge, so every original edge has cut probability
\[
 p=\frac13+\frac23\frac23=\frac79
\]
in an equal-part distribution on $C$.

Apply \cref{lem:flip-absorption}. Since the coefficient of $p_e$ in \eqref{eq:flip} is positive, every edge of $C$ has cut probability at least
\[
 \left(1-\frac4{5n}\right)\frac79+\frac2{5n}
 =\frac79-\frac2{9n}.
\]
For $n=6$ this is $20/27=\alpha+1/135$, and for $n=8$ it is $3/4=\alpha+1/60$. The claw edges have probability $\alpha$, proving $\U(C)$.
\end{proof}

In the matching contractions below, we retain all
nonmatching edge labels, obtaining a four-regular quotient
multigraph. Its underlying simple graph is the matching
quotient defined in \cref{sec:blocks}.

\begin{lemma}\label{lem:dummy}
Let $C$ be a simple triangle-free three-edge-connected cubic graph of order $n=2q$. There is an equal-part bisection distribution on $C\sqcup\claw$ cutting each claw edge with probability $\alpha$ and each edge of $C$ with probability at least
\[
 \beta_q=
 \begin{cases}
  \displaystyle\frac34-\frac3{20q},&q\text{ odd},\\[4pt]
  \displaystyle\frac34-\frac1{10q},&q\text{ even}.
 \end{cases}
\]
In particular $\U(C)$ holds if $q\geq9$ is odd or $q\geq6$ is even.
\end{lemma}
\begin{proof}
Choose a perfect matching $M$ of $C$ with edge marginals $1/3$, and contract its edges to obtain a loopless four-regular multigraph $Q$ on $q$ vertices, retaining all labels. For every nonempty proper set $S\subset V(Q)$, the cut $\delta_Q(S)$ is also a cut of $C$ and hence has size at least three. Since $Q$ is four-regular, this cut has even size, so $|\delta_Q(S)|\ge4$. Two edges between the same pair of contracted vertices have disjoint endpoints in $C$, by triangle-freeness, so their orientation constraints agree.

Set
\[
 \lambda_q=\begin{cases}
   \frac14,&q\text{ even},\\[2pt]
   \frac{q-1}{4q},&q\text{ odd}.
 \end{cases}
\]
The constant vector $x_h=\lambda_q$ on the labelled edges of $Q$ belongs to its matching polytope. Indeed, $x(\delta_Q(v))=4\lambda_q\le1$ and, for any proper odd set $S$,
\[
 x(E_Q(S))=\frac{\lambda_q}{2}\bigl(4|S|-|\delta_Q(S)|\bigr)
 \le 2\lambda_q(|S|-1)\le\frac{|S|-1}{2}.
\]
If $q$ is odd, the remaining odd set $V(Q)$ also satisfies the inequality, with equality $x(E(Q))=(q-1)/2$. Thus \cref{lem:multimatching} supplies a random matching $K$ of $Q$ with each label selected with probability $\lambda_q$. For even $q$, every realization of $K$ is perfect, since each quotient vertex has expected matching degree one.
For odd $q$, the expected size of $K$ is
$(q-1)/2$, the largest possible matching size.
Thus every realization leaves exactly one
quotient vertex $R$ unmatched. Moreover
\[
\Prb(R=v\mid M)=\frac1q
\qquad(v\in V(Q)).
\]

For each selected quotient edge, merge its two original matching pairs into one four-vertex block and orient it to cut all labels in that bundle. This is possible by the preceding compatibility observation. Orient distinct blocks by independent fair reversals; an unmatched quotient vertex remains a two-vertex block with its matched edge cut. This gives an equal-part bisection of $C$. Conditional on $M$, every matching edge is cut, while each nonmatching edge is cut with probability at least $1/2+\lambda_q/2$. Averaging over $M$, every edge of $C$ has cut probability at least
\[
 p_q=\frac13+\frac23\left(\frac12+\frac{\lambda_q}{2}\right)
 =\frac23+\frac{\lambda_q}{3}.
\]

If $q$ is even, $p_q=3/4$. Applying \cref{lem:flip-absorption} to this one distribution, with $n=2q$, yields an equal-part distribution on $C\sqcup\claw$ cutting each claw edge with probability $\alpha$ and each edge of $C$ with probability at least
\[
 \left(1-\frac2{5q}\right)\frac34+\frac1{5q}
 =\frac34-\frac1{10q}.
\]

Suppose $q$ is odd. We use the unmatched pair $R$ to perform the shift. Independently choose a balanced mode with probability $4/5$ and a shifted mode with probability $1/5$. In the balanced mode keep the bisection of $C$ and independently give the claw a uniformly random balanced cut. In the shifted mode, choose uniformly one of the two original vertices in the unmatched pair $R$ and flip its color; put the claw center on the majority side of $C$ and its three leaves on the minority side. Each outcome is a bisection of $C\sqcup\claw$, and each claw edge is cut with probability $(4/5)(2/3)+1/5=\alpha$.

Conditional on $M$, a matching edge becomes uncut
only when its pair is $R$ in the shifted mode.
Its cut probability is therefore $1-1/(5q)$.
 A nonmatching edge still has cut probability at least $1/2+\lambda_q/2$ in each mode. To see this in the shifted mode, observe that $R$ is not in a selected quotient bundle; the flipped pair becomes monochromatic with an unbiased color independent of the other fairly reversed blocks. Selected bundles remain fully cut, and endpoints in distinct blocks are separated with probability $1/2$. Consequently every edge of $C$ has cut probability at least
\[
 \frac13\left(1-\frac1{5q}\right)
 +\frac23\left(\frac12+\frac{q-1}{8q}\right)
 =\frac34-\frac3{20q}.
\]
The stated thresholds follow directly.
\end{proof}

Together, \cref{lem:dummy,lem:six-eight} cover all
even orders at least six except ten and fourteen.

\subsection{The ten-vertex case}

For ten-vertex graphs other than the Petersen graph,
compatible partitions of the five-vertex matching quotient
provide the required distribution.

\begin{lemma}\label{lem:five-quotient}
Let $C$ be a simple triangle-free cubic graph on ten vertices, let $M$ be a perfect matching, and contract the five edges of $M$ to obtain a loopless four-regular multigraph $Q$ on five vertices. If $C$ is not the Petersen graph, then $V(Q)$ has a finite distribution on partitions into blocks compatible with $(C,M)$ in which the endpoints of every quotient edge lie in one block with probability at least $1/3$.
\end{lemma}
\begin{proof}
Let $S$ be the simple underlying graph of $Q$,
let $D$ be the graph of double bundles of $Q$,
and put $H=K_5-S$. By triangle-freeness,
each bundle of $Q$ has size at most two.
Since $Q$ is four-regular,
\begin{equation}\label{eq:five-degree}
 d_D(v)=4-d_S(v)=d_H(v)\qquad(v\in V(Q)).
\end{equation}
The graphs $D$ and $H$ are edge-disjoint. By \eqref{eq:five-degree}, $2d_H(v)=d_D(v)+d_H(v)\leq4$, so $d_H(v)\leq2$ for every $v$.

Write $d=|E(H)|=|E(D)|$. The possibilities are
\begin{equation}\label{eq:five-supports}
 H\cong 5K_1,\qquad 2K_2\sqcup K_1,\qquad P_3\sqcup K_2,\qquad C_5.
\end{equation}
Indeed, $d\leq5$. The case $d=1$ is impossible because two one-edge graphs with the same vertex degrees must use the same edge. If $d=2$, a $P_3$ would force its degree-two vertex in $D$ to use the same two incident edges, so $H\cong 2K_2\sqcup K_1$. If $d=3$, a triangle would force the same triangle in $D$, while for a $P_4$ with isolated vertex an internal vertex has only one available neighbor outside its two $H$-edges; hence $H\cong P_3\sqcup K_2$. If $d=4$, the graph $H\cup D$ has eight edges and every vertex has even degree by \eqref{eq:five-degree}; its two-edge complement in $K_5$ would also have all degrees even, which is impossible. Finally $d=5$ forces every vertex of $H$ to have degree two, so $H\cong C_5$.

Suppose first that $d>0$. In each case of \eqref{eq:five-supports} there is a distribution on partitions into tree blocks for which every edge of $S$ is internal with probability at least $1/3$.
In the displays below, strings of digits denote vertex
sets, and vertical bars separate blocks of a partition.
If $H\cong 2K_2\sqcup K_1$, label the missing edges $01,23$ and use the three equiprobable partitions
\[
 014\mid2\mid3,\qquad 013\mid24,\qquad 012\mid34.
\]
If $H\cong P_3\sqcup K_2$, label the missing edges $01,12,34$ and use
\[
 04\mid23\mid1,\qquad 03\mid24\mid1,\qquad 02\mid134
\]
with equal probabilities. In each display, every edge of $S$ is internal in
exactly one of the three partitions.
If $H\cong C_5$, then $S$ is also a $5$-cycle;
averaging the partition consisting of three consecutive vertices
and the remaining adjacent pair over the five cyclic translates
with respect to the cyclic order of $S$
makes every edge internal with probability $3/5$.
Every nontrivial block displayed here induces a tree in $S$, and hence is compatible by \cref{lem:compatible-blocks}.

It remains to consider $S=K_5$.
In this case, compatibility also depends on
how the edges meet the original matching pairs.
Order the endpoints of the matching pair indexed by $i$ as $i_0,i_1$. For $i\ne j$, write the unique edge between the $i$- and $j$-pairs as
\[
 i_{a_{ij}}j_{a_{ji}},\qquad a_{ij}\in\{0,1\}.
\]
Here $\oplus$ denotes addition modulo $2$.
If the orientation variable $x_i\in\{0,1\}$ assigns
$i_t$ to side $t\oplus x_i$, then this edge is cut exactly when
\begin{equation}\label{eq:signed-constraint}
 x_i\oplus x_j=c_{ij},\qquad
 c_{ij}:=1\oplus a_{ij}\oplus a_{ji}.
\end{equation}
Reversing the order of one matching pair switches all constraints incident with the corresponding quotient vertex. Switch so that $c_{0i}=0$ for $1\leq i\leq4$, and let $R$ be the graph on $\{1,2,3,4\}$ whose edges are the pairs with $c_{ij}=1$.

Each vertex of $C$ has two neighbors outside $M$, so every row $(a_{ij})_{j\ne i}$ contains exactly two ones. Consequently, working modulo two and using $c_{0i}=0$,
\[
 \sum_{1\leq i<j\leq4}c_{ij}
 \equiv \sum_{i=1}^4\sum_{\substack{1\leq j\leq4\\j\ne i}}a_{ij}
 \equiv \sum_{i=1}^4 a_{i0}
 \equiv \sum_{i=1}^4(1+a_{0i})
 \equiv0.
\]
Thus $|E(R)|$ is even.

With $c_{0i}=0$, a block $B$ is compatible exactly when $R[B-\{0\}]$ is empty if $0\in B$, and when $R[B]$ is complete bipartite if $0\notin B$ (an empty part is allowed). This is the solvability criterion for \eqref{eq:signed-constraint}.

If $E(R)=\varnothing$, all five vertices form one compatible block.
Suppose $|E(R)|\in\{2,4\}$. Switching at a quotient vertex
complements its four incident signs and preserves block compatibility
(by changing the corresponding variable in \eqref{eq:signed-constraint}).
If $|E(R)|=4$, switch a degree-three vertex if one exists;
otherwise $R=C_4$, and after labeling its vertices
$1,2,3,4$ in cyclic order, switch vertices $2$ and $4$.
In the latter case, the four cycle signs become zero,
while $c_{02}=c_{04}=1$.
In either case the resulting signing has exactly
two pairs $ij$ with $c_{ij}=1$.
Relabel the quotient vertices so that neither pair
is incident with $0$.

If these two pairs are adjacent, label them $12,23$ and use
the three equiprobable partitions
\[
  0\mid123\mid4,\qquad 024\mid1\mid3,\qquad 0134\mid2.
\]
If they are disjoint, label them $12,34$. Use $0\mid12\mid34$
with probability $1/3$, and each partition $0xy\mid x'y'$,
where $x\in\{1,2\}$, $y\in\{3,4\}$ and $x',y'$ are their respective
partners, with probability $1/6$. Each displayed block is compatible,
and every pair occurs in one block with probability at least $1/3$.
Undoing the switches gives the desired distribution for the original
signing.

The only remaining possibility is $R=K_4$. We show that it forces $C$ to be the Petersen graph. Since the zeroth row has two ones, relabel $1,2,3,4$ so that
\[
 a_{01}=a_{02}=0,
 \qquad a_{03}=a_{04}=1.
\]
The normalization $c_{0i}=0$ gives $a_{i0}=1-a_{0i}$. Since $c_{ij}=1$ for every $1\leq i<j\leq4$, \eqref{eq:signed-constraint} gives $a_{ij}=a_{ji}$ for $1\leq i<j\leq4$. Let $L$ be the graph on $\{1,2,3,4\}$ in which $ij$ is an edge when $a_{ij}=1$. The row sums give
\[
 (d_L(1),d_L(2),d_L(3),d_L(4))=(1,1,2,2),
\]
so $L$ is a $P_4$ with endpoints $1,2$. The graph $C-M$ therefore consists of two $5$-cycles: one is obtained from $L$ by adjoining $0$ to its endpoints, and the other from the complementary $P_4$ by adjoining $0$ to its endpoints. Their edge sets are complementary $5$-cycles on the same label set, and $M$ joins equal labels. Thus $C$ is the Petersen graph. Hence $R=K_4$ cannot occur when $C$ is not the Petersen graph, and the lemma follows.
\end{proof}

\begin{lemma}\label{lem:ten-atom}
Every simple triangle-free three-edge-connected
cubic graph $C$ on ten vertices satisfies $\U(C)$.
More precisely, $C\sqcup\claw$ admits a finite
equal-part bisection distribution cutting every
edge with probability at least $\alpha$, whose
restriction to $C$ has part sizes $5$--$5$
with probability $4/5$ and $4$--$6$
with probability $1/5$.
\end{lemma}

\begin{proof}
Assume first that $C$ is not the Petersen graph. Choose a perfect matching $M$ from \cref{lem:multimatching}, so every edge belongs to $M$ with probability $1/3$. Conditional on $M$, apply \cref{lem:five-quotient,lem:compatible-blocks}. Every edge of $M$ is cut, while every other edge is cut with conditional probability at least
\[
 \frac12+\frac12\cdot\frac13=\frac23.
\]
Thus $C$ has a distribution on   $5$--$5$ bisections in which every edge is cut with probability at least
\[
 \frac13+\frac23\cdot\frac23=\frac79.
\]
Apply \cref{lem:flip-absorption} with $n=10$. It gives $\U(C)$ and the stated $5$--$5$ and $4$--$6$ probabilities on restriction to $C$, with every edge of $C$ cut with probability at least
\[
 \frac{23}{25}\cdot\frac79+\frac1{25}
 =\frac{34}{45}=\alpha+\frac1{45}.
\]

For the Petersen graph use the two-subset model on $\{1,2,3,4,5\}$, with adjacency defined by disjointness. The independent set $I=\{12,13,14,15\}$
has cut size twelve, while
$J=I\cup\{23\}$ has exactly two
internal edges and hence cut size eleven. The permutations of the five ground-set symbols act transitively on the fifteen edges. Therefore the uniformly permuted cuts $J$ and $I$ have edge marginals $11/15$ and $12/15$, respectively. Use the first law with probability $4/5$ and the second with probability $1/5$; every edge then has marginal
\[
 \frac45\cdot\frac{11}{15}+\frac15\cdot\frac{12}{15}
 =\frac{56}{75}>\alpha.
\]
The first law gives \(5\)--\(5\) bisections of \(C\),
whereas the second gives \(4\)--\(6\) partitions.

For this Petersen construction, in the $5$--$5$ mode give the extra claw a uniformly random balanced cut. In the $4$--$6$ mode cut all three claw arms and orient its three-vertex side towards the four-vertex side of $C$. Each union partition has seven vertices on each side, and each claw edge has probability $(4/5)(2/3)+1/5=\alpha$. This proves $\U(C)$.
\end{proof}

\subsection{The fourteen-vertex case}

\begin{lemma}\label{lem:seven-quotient}
Let $Q$ be a loopless four-regular multigraph on seven vertices, with multiplicity at most two between any pair of vertices, and let $S$ be its simple underlying graph. There is a finite distribution on partitions of $V(S)$ into singletons, edges, and induced three-vertex paths, in which the endpoints of every edge of $S$ lie in one block with probability at least $1/4$.
\end{lemma}
\begin{proof}
Let $d$ be the number of double bundles. Since $Q$ has fourteen labelled edges, $S$ has $14-d$ edges. For every five-vertex set $U$ with complementary pair $\{x,y\}$,
\begin{equation}\label{eq:seven-five}
 |E_Q(U)|=14-4-4+m_Q(xy)\leq8.
\end{equation}

\smallskip\noindent\emph{Case 1: at least two double bundles.}

If $d\geq2$, the constant vector $x_e=1/4$ belongs to the matching polytope of $S$.
 Vertex sums are at most one. A three-set spans at most three edges.
 For a five-set $U$, \eqref{eq:seven-five} gives
$|E_S(U)|\leq |E_Q(U)|\leq8$,
so its internal weight is at most two.
 The full seven-set has weight at most $12/4=3$, and singletons have no internal edges. This checks every odd-set inequality. Use the resulting matching distribution, taking its edges as two-vertex blocks and unmatched vertices as singletons.

\smallskip\noindent\emph{Case 2: no double bundle.}

If $d=0$, then $S=Q$ is simple four-regular and its complement $H$ is two-regular on seven vertices. It is therefore either a seven-cycle or a disjoint union of a three-cycle and a four-cycle.

If $H$ is a seven-cycle, label its vertices $0,\ldots,6$ cyclically. Use the partition
\[
 \{0,1,3\},\quad\{2,5\},\quad\{4,6\},
\]
and choose uniformly among its seven cyclic translates. The first block is an induced path in $S$, centered at $3$. Of the four internal edges, two have cyclic distance two and two have cyclic distance three. There are seven $S$ edges of each distance, and cyclic translation is transitive within each class. Thus each $S$ edge is internal with probability $2/7\geq1/4$.

If $H$ is a three-cycle on $A=\{a_1,a_2,a_3\}$ and a four-cycle on $B=\{b_1,b_2,b_3,b_4\}$, then $S$ contains all twelve $A$--$B$ edges, no edges within $A$, and two independent edges within $B$. Label the latter $b_1b_2,b_3b_4$. Start with
\[
 \{a_1,a_2,b_1\},\quad\{a_3,b_2\},\quad\{b_3,b_4\}.
\]
Average over all permutations of $A$ and all permutations of $B$ preserving its two-edge matching. The first block is an induced path. Each outcome makes three $A$--$B$ edges and one edge within $B$ internal. The permutation group is transitive on the twelve cross edges and on the two edges within $B$, so their internal probabilities are $1/4$ and $1/2$, respectively.

\smallskip\noindent\emph{Case 3: exactly one double bundle.}

Let its endpoints be $a,b$. The graph $S$ has thirteen edges; $a,b$ have degree three, the remaining vertices degree four, and $ab\in E(S)$.
We first reduce the argument to finding a spanning
forest $F$ consisting of a three-vertex induced path
in $S$ and two disjoint edges, with the path center
different from $a,b$.

Such a forest has four edges, one degree-two vertex, and six degree-one vertices. Put $R=S-E(F)$; it has nine edges and maximum degree at most three.
We claim that the constant vector $x_e=1/3$ belongs to the matching polytope of $R$.
 Degrees, three-sets and the full seven-set satisfy the required inequalities. Consider a five-set $U$ with complementary pair $\{x,y\}$. We have
\[
 |E_S(U)|=13-d_S(x)-d_S(y)+\mathbf1_{xy\in E(S)}.
\]
If $\{x,y\}=\{a,b\}$, this value is eight. Neither vertex is the center of $F$, so at most two forest edges meet the complementary pair; at least two forest edges are internal to $U$. Hence $|E_R(U)|\leq6$. If exactly one of $x,y$ is in $\{a,b\}$, then $|E_S(U)|\leq7$, and at most three forest edges meet the complementary pair. At least one forest edge is internal, again giving at most six remaining edges. If neither belongs to $\{a,b\}$, then $|E_S(U)|\leq6$ directly. These cases prove all five-set inequalities, completing the matching-polytope check.

With probability $1/4$ use the three components of $F$ as blocks; with probability $3/4$ use the matching distribution of $R$ and its unmatched singletons. An edge of $F$ is internal with probability $1/4$; an edge of $R$ is internal with probability $(3/4)(1/3)=1/4$. The induced-path condition ensures that no additional edge lies inside its three-vertex block.

We now construct $F$. Let $H=\overline S$. Then
$d_H(a)=d_H(b)=3$, every other vertex has degree two in $H$,
and $ab\notin E(H)$. First,
$N_H(a)\ne N_H(b)$: if these two neighborhoods were the same
three-element set, each of its vertices would already have both
its $H$-neighbors, leaving the other two vertices with degree at
most one in $H$.

Choose $c\in N_H(a)\setminus N_H(b)$, and write
$N_H(c)=\{a,t\}$. Choose any
$d\in N_H(b)\setminus\{t\}$. In $S$ we have
$bc,cd\in E(S)$ and $bd\notin E(S)$, so $b-c-d$ is an
induced path centered at $c\notin\{a,b\}$.

Put $T=S-\{b,c,d\}$. The graph $T$ has four vertices and no
isolated vertex. Indeed, $a$ has degree three in $S$, with
$ab\in E(S)$ and $ac\notin E(S)$, so it retains at least one
neighbor. Every other vertex of $T$ has degree four in $S$
and also retains at least one neighbor. Moreover, the two vertices
of $N_H(b)\setminus\{d\}$ belong to $T$ and each has degree
at least two in $T$, because neither is adjacent to $b$ in $S$.

The only four-vertex graph without isolated
vertices and without a perfect matching is
$K_{1,3}$. Since $T$ has at least two vertices
of degree at least two, it has a perfect matching.
 Together with the induced path $b-c-d$, its two
matching edges give the required forest $F$.
\end{proof}

\begin{lemma}\label{lem:fourteen-atom}
Every simple bridgeless triangle-free cubic graph $C$ on fourteen vertices has an equal-part bisection distribution cutting every edge with probability at least $3/4$. Consequently $\U(C)$ holds when $C$ is connected.
\end{lemma}
\begin{proof}
Choose a perfect matching $M$ with marginals $1/3$ and contract its seven pairs, retaining all nonmatching labels. The quotient is loopless four-regular. Between two pairs, triangle-freeness forces the cross edges to form a matching, so their multiplicity is at most two and their relative-orientation requirements agree. Apply \cref{lem:seven-quotient,lem:compatible-blocks}. Every matching edge is cut, and every nonmatching edge has conditional cut probability at least $1/2+(1/2)(1/4)=5/8$. Thus each edge of $C$ has probability at least
\[
 \frac13+\frac23\frac58=\frac34
\]
in one finite distribution on   $7$--$7$ bisections.

Apply \cref{lem:flip-absorption} with $n=14$. Its equal-part distribution on $C\sqcup\claw$ cuts each claw edge with probability $\alpha$ and each edge of $C$ with probability at least
\[
 \frac{33}{35}\cdot\frac34+\frac1{35}
 =\frac{103}{140}=\alpha+\frac1{420}.
\]
This proves $\U(C)$.
\end{proof}

Now, we are ready to give a proof of \cref{thm:claw-absorption}.

\begin{proof}[Proof of \cref{thm:claw-absorption}]
We first establish $\U(K)$ for every simple triangle-free three-edge-connected cubic graph $K$, using only the core constructions above. The cases of orders six and eight follow from \cref{lem:six-eight}, order ten from \cref{lem:ten-atom}, and order fourteen from \cref{lem:fourteen-atom}. The even-$q$ case of \cref{lem:dummy} covers orders divisible by four and at least twelve; its odd-$q$ case covers orders congruent to two modulo four and at least eighteen. Together these cover every possible even order at least six. None of these core lemmas uses the present absorption theorem.

Now let $C$ be any connected simple triangle-free cubic graph. If $C$ has a bridge, apply \cref{lem:bridged-absorption}. Otherwise, if it has a two-edge cut, apply \cref{lem:two-edge-reduction} with the established core assertion. If neither case applies, $C$ is three-edge-connected and is covered directly by the core assertion.
\end{proof}

\section{Parallel edges and sharpness}
\label{sec:components}

\Cref{thm:target} extends to finite loopless triangle-free subcubic multigraphs, with degrees counted with multiplicity. Replace every parallel bundle by a single edge, retaining the complete vertex set. The resulting simple graph is triangle-free and subcubic, and every vertex partition cuts either all labels of a bundle or none of them. The common distribution obtained in the proof of \cref{thm:target} therefore lifts to every labelled edge. In the weighted formulation, give the simple edge the sum of its bundle weights.

The only possible exceptional simple underlying graph is the entire claw. Its center already has three distinct neighbors, so the multiplicity-counting degree bound forces every arm to have multiplicity one. The original multigraph would itself be the claw. Thus no new exception is introduced.

\begin{proposition}\label{prop:sharp}
The Petersen graph has fifteen edges and maximum bisection
size eleven. The claw has maximum bisection size two.
\end{proposition}
\begin{proof}
Represent the Petersen graph by the two-element subsets of $\{1,2,3,4,5\}$, with adjacency defined by disjointness. It has ten vertices, degree three, and fifteen edges. It has no triangle because there are no three pairwise disjoint two-element subsets of a five-element set.

Any four distinct pairwise intersecting two-element sets have a common element. Indeed, if there is no common element, choose intersecting sets $\{a,b\},\{a,c\}$ and a set avoiding $a$. That set must be $\{b,c\}$, and every set meeting all three is one of those three. Thus an independent four-set of Petersen vertices consists of the four pairs containing one fixed element.

Every five-vertex set $S$ in the Petersen graph spans at least two edges. If it spanned at most one, deleting an endpoint of that edge (or any vertex if there were no edge) would leave an independent four-set, the four pairs through one fixed element. The fifth pair omits that element and is disjoint from exactly two of those four pairs. This forces at least two edges in $S$, a contradiction. Therefore every bisection has at most
\[
 |\delta(S)|=3|S|-2|E(S)|\leq15-4=11
\]
edges. Taking the four pairs through a fixed element and any fifth pair omitting it gives exactly two internal edges, and therefore a bisection of size eleven.

In a $2$--$2$ partition of the claw, its center shares its part with exactly one leaf, so exactly two of its three arms are cut. The universal constant cannot exceed $11/15$ on the target class, and excluding the entire claw is necessary.
\end{proof}

\section*{Declaration on the Use of Generative AI}
The authors used generative AI tools to assist in discussing proof strategies, checking proofs, and improving
exposition.

\end{document}